\documentclass[12pt,leqno]{article}

\usepackage{graphicx}
\font\small=cmr10

\def\Section{\setcounter{equation}{0} \setcounter{thm}{0} 
\setcounter{lem}{0} \setcounter{cor}{0} \setcounter{defcor}{0} 
\setcounter{defi}{0} \setcounter{prop}{0} \setcounter{remark}{0} 
\setcounter{example}{0} \setcounter{prob}{0} \section}
\def\thesection       {\arabic{section}.}

\newlength{\noulong}
\newlength{\noulonga}
\font\itpetit=cmti10

\def\beq{\begin{equation}\displaystyle}
\def\eeq{\end{equation}}
\def\bel{\begin{equation} \displaystyle \begin{array}{l} }
\def\eel{\end{array} \end{equation} }
\def\bell{\begin{equation} \displaystyle \begin{array}{ll}  }
\def\eell{\end{array} \end{equation} }

\def\bea{\begin{eqnarray}}
\def\eea{\end{eqnarray} }
\def\bean{\begin{eqnarray*}}
\def\eean{\end{eqnarray*} }

\def\bar#1{\overline#1}

\def\T{{\bf T}}

\def\TT{{\cal T}}
\def\CC{\rm \hbox{C\kern-.57em\raise.47ex
         \hbox{$\scriptscriptstyle |$}\kern+0.5 em }}

\def\NN{{\rm I\hspace{-0.50ex}N} }
\def\OO{\rm \hbox{O\kern-.34em\raise.47ex
         \hbox{$\scriptscriptstyle |$}\kern-.46em\raise.47ex
         \hbox{$\scriptscriptstyle |$}\kern+0.5 em }}
\def\RR{{\rm I\hspace{-0.50ex}R} }

\def\debmin{
   \hspace{1mm} \smallskip\newline
   \hspace*{20mm}\begin{minipage}{\noulonga}}
\def\finmin{\end{minipage}}

\catcode`@=11
\def\eqalign#1{\null\,\vcenter{\openup1\jot \m@th
   \ialign{\strut \hfil$\displaystyle{##}$ & $\displaystyle{{}##}$\hfil
      \crcr#1\crcr}}\,}
\catcode`@=12

\newtheorem{thm}{Theorem}  
\newtheorem{lem}{\bf Lemma}
\newtheorem{cor}{\bf Corollary}
\newtheorem{defcor}{\bf Definition and Corollary}
\newtheorem{defi}{\bf Definition}
\newtheorem{prop}{\bf Proposition}
\newtheorem{remark}{\bf Remark} 
\newtheorem{example}{\bf Example}
\newtheorem{prob}{\bf Problem}

\def\augmentetout{\refstepcounter{thm}\refstepcounter{lem}\refstepcounter{cor}\refstepcounter{defcor}\refstepcounter{prob}\refstepcounter{prop}\refstepcounter{defi}\refstepcounter{remark}\refstepcounter{example}}

\def\finproof{\mbox{}\hfill {\small \fbox{}} \\}

\def\debthm#1#2finthm {\begin{thm} #1 #2 \augmentetout \addtocounter{thm}{-1}\end{thm} }

\def\deblem#1#2finlem {\begin{lem} #1 #2 \augmentetout \addtocounter{lem}{-1}\end{lem} }

\def\debprop#1#2finprop {\begin{prop}#1 #2 \augmentetout \addtocounter{prop}{-1}\end{prop} }

\def\debcor#1#2fincor {\begin{cor}#1 #2 \augmentetout \addtocounter{cor}{-1}\end{cor} }

\def\debdefcor#1#2findefcor {\begin{defcor}#1 #2 \augmentetout \addtocounter{defcor}{-1}\end{defcor} }

\def\debprob#1#2finprob{\begin{prob}#1 #2  \augmentetout \addtocounter{prob}{-1}\end{prob} }

\def\debdef#1#2findef {\begin{defi}#1 #2  \augmentetout \addtocounter{def}{-1}\end{defi} }

\def\debrem#1#2finrem{ \begin{remark}#1 #2 \augmentetout \addtocounter{rem}{-1}\end{remark}}

\def\debex{\begin{example} \rm\mbox{}\augmentetout \addtocounter{example}{-1}\\}
\def\finex{\end{example}}

\def\debbox{ 
   \hspace{10mm}\smallskip\newline 
   \hspace*{3mm}\begin{tabular}[b]{|r |l} \hline \\
   \rm\makebox[2mm]{}\begin{minipage}{\noulong}}
\def\finbox{\end{minipage}\\ \\ \hline \end{tabular}}

\catcode`@=11
\def\thebibliography#1{\noindent {\bf References}\sloppy\nolinebreak
 \list
 {[\arabic{enumi}]}{\settowidth\labelwidth{[#1]}\leftmargin\labelwidth
 \advance\leftmargin\labelsep
 \usecounter{enumi}}
 \def\newblock{\hskip .11em plus .33em minus -.07em}
 \sloppy
 \sfcode`\.=1000\relax}
\catcode`@=12
\catcode`@=11
\def\chaptermark#1{\markboth {\itpetit \ifnum \c@secnumdepth >\m@ne
 \@chapapp\ \thechapter. \ \fi #1}{}}
\def\sectionmark#1{\markright {\itpetit \ifnum \c@secnumdepth >\z@
 \thesection \ \ \fi  #1  }}
\catcode`@=12
\newcounter{clegende}
\newread\epsffilein    
\newif\ifepsffileok    
\newif\ifepsfbbfound   
\newif\ifepsfverbose   
\newdimen\epsfxsize    
\newdimen\epsfysize    
\newdimen\epsftsize    
\newdimen\epsfrsize    
\newdimen\epsftmp      
\newdimen\pspoints     
\def\epsfbox#1{\global\def\epsfllx{72}\global\def\epsflly{72}%
   \global\def\epsfurx{540}\global\def\epsfury{720}%
   \def\lbracket{[}\def\testit{#1}\ifx\testit\lbracket
   \let\next=\epsfgetlitbb\else\let\next=\epsfnormal\fi\next{#1}}%
\def\epsfgetlitbb#1#2 #3 #4 #5]#6{\epsfgrab #2 #3 #4 #5 .\\%
   \epsfsetgraph{#6}}%
\def\epsfnormal#1{\epsfgetbb{#1}\epsfsetgraph{#1}}%
\def\epsfgetbb#1{%
\openin\epsffilein=#1
\ifeof\epsffilein\errmessage{I couldn't open #1, will ignore it}\else
   {\epsffileoktrue \chardef\other=12
    \def\do##1{\catcode`##1=\other}\dospecials \catcode`\ =10
    \loop
       \read\epsffilein to \epsffileline
       \ifeof\epsffilein\epsffileokfalse\else
          \expandafter\epsfaux\epsffileline:. \\%
       \fi
   \ifepsffileok\repeat
   \ifepsfbbfound\else
    \ifepsfverbose\message{No bounding box comment in #1; using defaults}\fi\fi
   }\closein\epsffilein\fi}%
\def\epsfclipstring{}
\def\epsfsetgraph#1{%
   \epsfrsize=\epsfury\pspoints
   \advance\epsfrsize by-\epsflly\pspoints
   \epsftsize=\epsfurx\pspoints
   \advance\epsftsize by-\epsfllx\pspoints
   \epsfxsize\epsfsize\epsftsize\epsfrsize
   \ifnum\epsfxsize=0 \ifnum\epsfysize=0
      \epsfxsize=\epsftsize \epsfysize=\epsfrsize
      \epsfrsize=0pt
     \else\epsftmp=\epsftsize \divide\epsftmp\epsfrsize
       \epsfxsize=\epsfysize \multiply\epsfxsize\epsftmp
       \multiply\epsftmp\epsfrsize \advance\epsftsize-\epsftmp
       \epsftmp=\epsfysize
       \loop \advance\epsftsize\epsftsize \divide\epsftmp 2
       \ifnum\epsftmp>0
          \ifnum\epsftsize<\epsfrsize\else
             \advance\epsftsize-\epsfrsize \advance\epsfxsize\epsftmp \fi
       \repeat
       \epsfrsize=0pt
     \fi
   \else \ifnum\epsfysize=0
     \epsftmp=\epsfrsize \divide\epsftmp\epsftsize
     \epsfysize=\epsfxsize \multiply\epsfysize\epsftmp   
     \multiply\epsftmp\epsftsize \advance\epsfrsize-\epsftmp
     \epsftmp=\epsfxsize
     \loop \advance\epsfrsize\epsfrsize \divide\epsftmp 2
     \ifnum\epsftmp>0
        \ifnum\epsfrsize<\epsftsize\else
           \advance\epsfrsize-\epsftsize \advance\epsfysize\epsftmp \fi
     \repeat
     \epsfrsize=0pt
    \else
     \epsfrsize=\epsfysize
    \fi
   \fi
   \ifepsfverbose\message{#1: width=\the\epsfxsize, height=\the\epsfysize}\fi
   \epsftmp=10\epsfxsize \divide\epsftmp\pspoints
   \vbox to\epsfysize{\vfil\hbox to\epsfxsize{%
      \ifnum\epsfrsize=0\relax
        \includegraphics{#1}%
      \else
        \epsfrsize=10\epsfysize \divide\epsfrsize\pspoints
        \includegraphics{#1}%
      \fi
      \hfil}}%
\global\epsfxsize=0pt\global\epsfysize=0pt}%
{\catcode`\%=12 \global\let\epsfpercent=
\long\def\epsfaux#1#2:#3\\{\ifx#1\epsfpercent
   \def\testit{#2}\ifx\testit\epsfbblit
      \epsfgrab #3 . . . \\%
      \epsffileokfalse
      \global\epsfbbfoundtrue
   \fi\else\ifx#1\par\else\epsffileokfalse\fi\fi}%
\def\epsfempty{}%
\def\epsfgrab #1 #2 #3 #4 #5\\{%
\global\def\epsfllx{#1}\ifx\epsfllx\epsfempty
      \epsfgrab #2 #3 #4 #5 .\\\else
   \global\def\epsflly{#2}%
   \global\def\epsfurx{#3}\global\def\epsfury{#4}\fi}%
\def\epsfsize#1#2{\epsfxsize}
\def\T{\mathbf{T}}
\def\beq{\begin{equation}}
\def\eeq{\end{equation}}

\def\di{\displaystyle}

\def\R2+{\RR ^2_+}

\def\ni{\noindent}

\def\M{{\cal M}}
\def\D{{\cal D}}

\def\NN{I\!\!N}
\def\RR{I\!\!R}

\def\intR2+{\int_{\R2+}}

\def\TT{{\cal T}}

\def\lto{\ensuremath{\longrightarrow}}
\def\A{{\cal A}}

\def\1{1\!\!1}

\begin{document}
\begin{center}

\centerline{\Large \bf Convergence to time-periodic solutions in time-periodic}
\centerline{\Large \bf Hamilton-Jacobi equations on the circle}

\bigskip

{\bf Patrick BERNARD}

\medskip {\footnotesize 
Institut Fourier, UMR CNRS 5582
Universit\'e de Grenoble I 
BP 74,\\ 38402 Saint-Martin d'H\`eres, France 
}

\bigskip

{\bf Jean-Michel ROQUEJOFFRE}

\medskip {\footnotesize UFR-MIG, UMR CNRS 5640.  Universit\'e Paul
Sabatier, 118 route de Narbonne,\\ 31062 Toulouse Cedex, France}

\rule{30mm}{.1mm}
\end{center}

\ni {\small {\sc Abstract.} \rule [.5mm]{2mm}{.1mm} 
.} {\small The goal of this paper is to give a simple  proof of the convergence to
time-periodic states  of the solutions of time-periodic
Hamilton-Jacobi equations on the circle with convex Hamiltonian.
Note that the period of limiting solutions may be greater than the period 
of the Hamiltonian.}

\ni {\small {\sc Resum\'e.} \rule [.5mm]{2mm}{.1mm} 
.} {\small On donne une preuve simple de la convergence vers des \'etats
p\'eriodiques en temps pour les solutions d'\'equations de Hamilton-Jacobi
sur le cercle avec un  Hamiltonien convexe et p\'eriodique en temps.
Il est \`a noter que la p\'eriodes des solutions limites
peut \^etre plus grande que la p\'eriode du Hamiltonien.}

\Section{Introduction}
We consider the Hamilton-Jacobi equation
\beq
\label{HJ}
u_t+H(t,x,u_x)=0,\ \ \ x\in \T
\eeq
where $\T$ is the unit circle. 
The Hamiltonian $H(t,x,p):\ \RR\times\T\times\RR\mapsto\RR$
is $C^2$, 1-periodic in $t$, 
and satisfies the following classical hypotheses:
\begin{itemize}
\item
Strict convexity: $H_{pp}(t,x,p)>0$ for all 
$(t,x,p)\in\RR\times\T\times\RR.$
\item
Super-linearity: $H(t,x,p)/p\lto \infty$ as $|p|\lto \infty$
for each $(t,x)\in \RR\times\T$.
\item
Completeness: The Hamiltonian vector-field
$$X(t,x,p)=(H_p(t,x,p),-H_x(t,x,p))$$
is complete, \textit{i.e.}
for all $(t_0,x_0,p_0)$, there exists a $C^2$ curve
$\gamma(t)=(x(t),p(t)):\RR\lto \T\times \RR$
such that $(x(t_0),p(t_0))=(x_0,p_0)$
and $\dot \gamma(t)=X(t,\gamma(t))$ for all $t\in\RR$.
\end{itemize}
The first two assumptions are classical in the viscosity solutions theory;
see \cite{PL2}. The last one was introduced  in Mather \cite{Mather}; note that it is satisfied if
there exists a constant $C$ such that
$
|H_t|\leq C(1+H)$.

Under the above three assumptions, the Cauchy Problem for (\ref{HJ}) is well posed in the viscosity sense:
given a time $s\in \RR$ and a continuous function $u_0:\ \T\to\RR$, equation
(\ref{HJ}) has a unique viscosity solution 
$u(t,x):[s,+\infty[\times\T\lto \RR$, such that $u(s,.)=u_0$.
It will be denoted by
$\TT(s,t)u_0$. See \cite{PL2}, for instance.

It is known - and this is not specific to the one-dimensional setting - 
that 
there exists a real  number $\lambda$
such that 
$u(t,x)+\lambda t$
is bounded for all viscosity solution $u:[s,+\infty[\times\T\lto \RR$
of (\ref{HJ}).
The real number $\lambda$ has various  different names:
It is the critical value of Ma\~ne,
see \cite{Mane}, \cite{CDI} or the value 
$\alpha(0)$, where $\alpha:\RR\lto\RR$ is the 
Mather function, see \cite{Mather}, or the averaged Hamiltonian \cite{LPV}.
Note that the number $\lambda$ may also be viewed
as the eigenvalue of the Hopf-Lax-Oleinik operator
in the sense of min-plus algebra, see \cite{S}.
We are  interested in proving the following result.

\begin{thm}
\label{t1.1}
Let $u(t,x):[s,+\infty[\times\T\lto \RR$ be a viscosity solution
of (\ref{HJ}).
There exist an integer $T$ and a viscosity solution 
$l(t,x)=\phi(t,x)-\lambda t:\RR\times\T\lto\RR$
such that $\phi$ is $T$-periodic in $t$ and
$$
\lim _{t\lto \infty }\big\|u(t,.)-l(t,.)\big\|_{\infty}= 0.
$$
\end{thm}

In the following, we will always assume
that $\lambda=0$, which can be obtained by replacing
the Hamiltonian $H$
by $H-\lambda$.

It is known that there always exists a viscosity solution
of (\ref{HJ}) which is $1$-periodic in time
(such a solution is not unique in general).
However, it is not hard to build examples of
viscosity solutions of equations of the form 
(\ref{HJ}) which do not converge to $1$-periodic
solutions, see  \cite{BS2} and \cite{FaMa}.
More precisely, one can build solutions
which are periodic in time, but 
of minimal period greater than one.
Hence one cannot expect to have 
always $T=1$ in the theorem.

For time-independent Hamiltonians, convergence to steady states is known: a particular nontrivial multidimensional case (including some non-strictly convex situations) is studied in 
\cite{NR2}, and the general result in one dimension for strictly convex Hamiltonians 
is given in \cite{R2}. Recall - \cite{BS2} - that non strict convexity may result in a failure of the convergence result. The general multidimensional result in the strictly convex case is due to Fathi \cite{Fa3}, with dynamical systems arguments. A purely PDE proof, encompassing both assumptions of \cite{Fa3} and \cite{NR2}, is provided in \cite{BS2}.
 A proof mixing PDE and dynamical systems arguments is given in \cite{R}, and proves the above-stated theorem in some cases (rational rotation number; see below).

The situation is not so clear when the Hamiltonian is time-periodic.
In order to be more precise, it is useful to
recall that one can associate to the equation
(\ref{HJ}) a rotation number $\rho\in \RR$, see section \ref{circle},
  which is the 
rotation number of extremals.
We then have the following

\textbf{Addendum}
\begin{itshape}
The period $T$ 
in the theorem is $1$ if the rotation number $\rho$
is irrational and is not greater  than $q$ if 
$\rho$ is a rational $p/q$.
\end{itshape}

The theorem and its addendum have been proved in 
 \cite{Be} from the dynamical system point of view,
that is from the study of extremals.
Let us mention also that, 
in the case of an irrational rotation number,
there is a single periodic viscosity solution,
(up to an additive constant)
as was proved in \cite{E} and \cite{S}.

As said before, the result in the case of a rational rotation number
had been previously obtained by a method
relying on the dynamic programming principle
in  \cite{R}.
It turns out that the ideas of this paper
can be exploited further to provide
a simpler proof of the general case.
It is our aim to present this proof here.

We shall first recall the general properties
of viscosity solutions, in Section \ref{general}.
In Section \ref{calibrated}, we introduce
some dynamics, define the Aubry-Mather sets,
and recall 
specific observations concerning $\omega$-limit solutions,
mostly taken from \cite{R}.
All the results in these sections are general, 
and remain true if one considers equation
(\ref{HJ}) on any compact manifold.
We complement these general observations
by specific one-dimensional arguments
in Section \ref{circle} to conclude the proof.

{\bf Acknowledgement.} This paper is the result of conversations 
initiated at the june 2002 workshop on Hamilton-Jacobi equations,
held in Cortona. 
It is our pleasure to thank the organizers: F. Camilli, I. Capuzzo Dolcetta and  A. Siconolfi.
\Section{General properties and a large time behaviour candidate}
\label{general}
Let us start with some well understood properties of viscosity solutions without proof, see \cite{PL2}, \cite{Fa1}
or \cite{FaMa}.
It is useful to introduce the Lagrangian
$$
L(t,x,v)=\max_{p\in\RR} \biggl(pv-H(t,x,p)\biggl),
$$
which is also  convex and superlinear.
We have the following fundamental result of the calculus of variations
(see a proof in the appendix of \cite{Mather}):

Let us fix a time interval $[t,t']$, and two points
$X,X'$ in $\T$.
If the Lagrangian is convex and superlinear,
the action integral
$$
\int _t ^{t'} L(s,x(s),\dot x(s)) ds
$$
reaches its minimum on the set of 
absolutely continuous curves $x(s):[t,t']\lto  \T$
which satisfy $x(t)=X$ and $x(t')=X'$.
If in addition the Hamiltonian flow  is complete, then the minimum
is reached by $C^2$ curves which  satisfy the Euler-Lagrange equations.
Denoting 
$
p(s)= L_v(q(s),\dot q(s),s),$
these equations can be written
$$
\dot p(s)= L_x(s,x(s),\dot x(s))=
 -H_x(s,x(s),p(s))
$$
hence the curve $(x(s),p(s))$ is a trajectory of 
the Hamiltonian vector-field.

Let us define, for all $t>s$, the 
Hopf-Lax-Oleinik operator $\TT(s,t)$, which, to any continuous 
function $u_0:\T\lto \RR$, associates the continuous function

\beq
\label{LO}
\TT(s,t)u_0(x)=\inf_{\gamma}\biggl(u_0(\gamma(s))+\int_s^tL(\sigma,\gamma,\dot\gamma)\ d\sigma\biggl),
\eeq
where the infimum is taken on the set of absolutely continuous 
curves $\gamma:[s,t]\lto \T$ satisfying $\gamma(t)=x$.

The viscosity solution  $u(t,x)\in C([s,+\infty[\times\T,\RR)$
of (\ref{HJ}) with initial condition 
 $u(s,.)=u_0$
is given by $u(t,.)=\TT(0,t)u_0$.
This may be taken as the definition of viscosity solutions
in the present paper. This definition 
is equivalent to the standard one in our setting.
This fact is classical, the reader may find a
proof in \cite{FaMa} for instance.

We have the Markov property
$$
 \TT(t,t')\circ\TT(s,t)= \TT(s,t')
$$
for $s\leq t\leq t'$,
hence the mappings $\TT(0,n)=\TT(0,1)^n$, $n\in\NN$
form a discrete semi-group.
We will note $\TT$ for $\TT(0,1)$ for simplicity.
The following simple properties are proved in \cite{Fabook}
in the autonomous case, but the extension to the nonautonomous
case does not present any difficulty.

The mappings $\TT(s,t)$ are contractions,
$$
\big\|
\TT(s,t)u-\TT(s,t)v
\big\|_{\infty}
\leq
\big\|u-v
\big\|_{\infty}.
$$

The mappings  $\TT(s,t)$ are compact.
More precisely, 
there exists a positive non-increasing function
$K(\epsilon):]0,\infty[\lto]0,\infty[$
such that $\TT(s,t)u$
is $K(\epsilon)$-Lipschitz for all continuous function $u$ 
and all $t\geq s+\epsilon$.

The mappings  $\TT(s,t)$ are order-preserving and satisfy
$\TT(s,t)(c+u)=c+\TT(s,t)(u)$ for all real $c$.

The following proposition will be useful in the sequel.
Note that it  implies the existence of a fixed point of
the semi-group generated by $\TT$,
\textit{i.e.} the existence of a $1$-periodic viscosity solution
of (\ref{HJ}).

\begin{prop}
\label{liminf}
Let $u(t,x)\in C([s,\infty[\times \T,\RR)$
be a viscosity solution of (\ref{HJ}); recall that it is bounded on 
$[s,\infty[\times\T$.
Then the $1$-periodic function
$$
\phi(t,x)=\liminf_{n\to+\infty}u(t+n,x)
$$
is a viscosity solution of (\ref{HJ}).
\end{prop}

From the Barles-Perthame lemma \cite{BP},
$\phi$ is a viscosity sub-solution. Much more, due to the convexity of the Hamiltonian, it is
a solution - Barron-Jensen \cite{BJ}. We thank the referee for pointing  out to us the last item.
For the sake of self-containedness, we provide a simple proof using
the Hopf-Lax-Oleinik expression  of solutions.

\noindent{\sc Proof.} 
We have to prove that $\TT(s,t)\phi(s,.)=\phi(t,.)$
for all $s\leq t$.
Let us first prove that $\phi$ is a sub-solution, \textit{i. e. }
that  $\TT(s,t)\phi(s,.)\leq \phi(t,.)$.
In order to do so, we 
fix $(t,x)$ and consider 
 an increasing sequence $n_k$
of integers such that $u(t+n_k,x)\lto \phi(t,x)$.
There exists a sequence of curves $\gamma_k:[s,t]\lto\T$ such that
$$
u(t+n_k,x)=u(s+n_k,\gamma_k(0))+
\int_s^{t}L(\sigma,\gamma_k(\sigma+t),\dot\gamma_k(\sigma+t))\ d\sigma.
$$
The sequence $\gamma_k$ is compact for the $C^1$ topology,
and we will assume by possibly taking a subsequence in $n_k$
that it is convergent, and note $\gamma$ the limit.
Taking the $\liminf$ in the equality above gives 
$$
\phi(t,x)\geq \phi(s,\gamma(0))+
\int_s^{t}L(\sigma,\gamma(\sigma+t),\dot\gamma(\sigma+t))\ d\sigma
\geq \TT(s,t)\phi(s,.)(x).
$$
We have used that the functions $u(t,.)$ are equicontinuous
to conclude that $\liminf u(s+n_k, \gamma_k(0))=
\liminf u(s+n_k, \gamma(0))\geq \phi(s,\gamma(0))$.

We now prove that  $\phi$ is a super-solution, \textit{i. e. }
that  $\TT(s,t)\phi(s,.)\geq \phi(t,.)$.
Note that for all curve $\gamma:[s,t]\lto\T$,
we have 
$$
u(t+n,x)\leq u(s+n,\gamma(0))+
\int_s^{t}L(\sigma,\gamma(\sigma+t),\dot\gamma(\sigma+t))\ d\sigma.
$$
Taking the $\liminf$, we obtain 
$$
\phi(t,x)\leq \phi(s,\gamma(0))+
\int_s^{t}L(\sigma,\gamma(\sigma+t),\dot\gamma(\sigma+t))\ d\sigma
$$
for each curve $\gamma$,
hence $\phi(t,.)\leq \TT(s,t)\phi(s,.)$,
which is the desired inequality.
 \finproof

The basic objects to understand in order to
study  the asymptotic behaviour of
solutions of  (\ref{HJ}) are the $\omega$-limit solutions.
Recall that a solution $u(t,x):\RR\times \T\lto \RR$
is called an $\omega$-limit
solution if there exists a solution 
$v:[s,\infty[\times \T\lto\RR$
and an increasing sequence $n_k$ of integers such that
$$u(t,x)=\lim_{k\lto \infty}v(t+n_k,x).$$
In other words, $\omega$-limit solutions are solutions
whose initial value $u(0,.)$ is an $\omega$-limit
of the semi-group $\TT$.
If $u(t,x)\RR\times \T\lto \RR$ 
is an $\omega$-limit solution of  (\ref{HJ}),
then it is bounded, and the functions
$u(t,.),t\in \RR$ are equilipschitz.


%
%
%
%
%
%
%

\section{Calibrated curves and Uniqueness set}
\label{calibrated}
In this section, we give   
some salient features of $\omega$-limit solutions that do not depend of the dimension of the ambient space.
The key differentiability result  is stated and used as a black box; the other results are proved, although their proofs are already in \cite{Fa1} or \cite{R}.

\medskip
Let $u:[s,\infty[\times\T\lto\RR$
be a viscosity solution of (\ref{HJ}).
A curve $\gamma:[s,\infty[\supset[t,t']\lto \T$
is said calibrated by $u$ if
$$
u(t',\gamma(t'))=u(t,\gamma(t))+
\int_{t}^{t'}L(\sigma,\gamma,\dot\gamma)\ d\sigma.
$$ 
Note that if $\gamma(s)$ is a calibrated curve,
then the curve $(\gamma(s),p(s))$ is a trajectory
of the Hamiltonian flow, where 
$p(s)=\partial_v L(s,\gamma(s),\dot \gamma(s))$.
We have the following regularity result
\begin{thm}
\label{threg} (Fathi,  \cite{Fa1}).
If $\gamma:[t,t']\lto \T$
is calibrated by $u$,
then $u_x$ exists at each point $(s,\gamma(s)),s\in [t,t'[$
and satisfies
\begin{equation}
\label{differentielles}
\begin{array}{rll}
&u_x(s,\gamma(s))=p(s)=L_v(s,\gamma(s),\dot\gamma(s))\\
\Longleftrightarrow&\\
&\dot\gamma(s)=H_p(s,\gamma(s),u_x(s,\gamma(s)).
\end{array}
\end{equation}
\end{thm}

It was already proved in \cite{Fleming}
that this relation holds at differentiability points of $u$.
The regularity of $u$  on calibrated curves 
was noticed  in \cite{PR} for the equation
$\vert\nabla u\vert=f(x)$ on the sphere, 
$f$ nonnegative with nonempty zero set.
It was not, however, made as systematic as in \cite{Fa1}.

\medskip
We now choose once and for all a 
 1-periodic solution  $\phi(t,x)$ of (\ref{HJ}).
First, let us note that a classical compactness argument gives the existence of
curves $\gamma:\RR\lto\T$
which are calibrated by $\phi$ on all compact interval.

As a consequence of Theorem \ref{threg},
two such curves cannot intersect.
Indeed, if $\gamma_1$ and  $\gamma_2:\RR\lto \T$
are calibrated by $\phi$, and if there exists a $t$ such that
$\gamma_1(t)=\gamma_2(t)$, 
then $u_x$ exists at the point $(t,\gamma_1(t))$,
and, setting 
$p_i(s)=\partial L(s,\gamma_i(s),\dot \gamma_i(s))$,
we have 
$$
p_1(t)=u_x(t,\gamma_1(t))=u_x(t,\gamma_2(t))=p_2(t).
$$
The curves $(\gamma_i(s),p_i(s))$ are then two trajectories
of the Hamiltonian flow which are equal at time $t$, 
so that they coincide for all times.

Let  $$ \D\subset \RR\times \T$$
be the union of the graphs of these orbits,
and $\D_0\subset \T$
be the set of points $\gamma(0)$, where 
$\gamma:\RR\lto \T$ is calibrated.
This is a nonempty compact set.
For each $t$, we define the mapping $S^t:\D_0\lto \T$
which associates to each $x\in \D_0$,
the value $\gamma(t)$, where $\gamma:\RR\lto \T$
is the unique calibrated curve satisfying $\gamma(0)=x$.
It is a  homeomorphism onto its image.
Let us mention for completeness that this 
homeomorphism is bi-Lipschitz, as can be obtained
from refined versions of Theorem \ref{threg}.
The first results in that direction are due to Mather.

Clearly, $S^1$ is a homeomorphism of $\D_0$.
Let us note $\M_0$ its $\omega$-limit.
This is the closure in $\T$ 
of the set of points $x\in \D_0$
which are the limit of a sequence 
$S^{n_k}(y)$ with $y\in \D_0$ and 
$n_k$ an increasing sequence of integers.
The set $\M_0$ is non-empty and compact.
We call $\M$ the union, in $\RR\times\T$,
of the graphs  of curves $S^t(x),x\in \M_0$.

\subsection{The behaviour of global solutions on global extremals}
Recall that we have fixed a periodic viscosity solution $\phi$.
The following remark, noticed in 
\cite{R}, is the key point to the convergence proof:

\begin{lem}
 \label{decroit}
Let $u_1$ and  $u_2(t,x):[s,\infty[\times \T\lto \RR$
be viscosity solutions  of (\ref{HJ}), and let $\gamma(t):[s,s']\lto \T$
be a curve calibrated by $u_2$,
Then
the function
 $ t\mapsto u_1(t,\gamma(t))-u_2(t,\gamma(t))$
 is non-increasing.
\end{lem}
\noindent{\sc Proof. } For every $t\leq t'$ we have:
$$
u_1(t',\gamma(t'))-u_1(t,\gamma(t))
\leq\int_t^{t'}L(s,\gamma(s),\dot \gamma(s))\ ds
$$ 
by the Lax-Oleinik formula, and
$$
u_2(t',\gamma(t'))-u_2(t,\gamma(t))=\int_t^{t'}L(s,\gamma(s),\dot \gamma(s))\ ds.
$$
Hence we have $(u_1-u_2)(t',\gamma(t'))\leq(u_1-u_2)(t,\gamma(t))$. 
\finproof
The important consequence below is also proved in \cite{R}
(recall that $\phi$ is a prescribes periodic viscosity solution):
\begin{cor}
\label{const}
Let $u(t,x)$ be an $\omega$-limit viscosity solution of  (\ref{HJ}),
and let $x\in \M_0$,
then 
the function  $ t\mapsto u(t,S^t(x))-\phi(t,S^t(x))$
is constant.
\end{cor}
\noindent {\sc Proof.} 
Choose $t_0>0$. Let $u_0\in C(\T)$ be such that $u$ is an $\omega$-limit solution to $u_0$.
Consider $x_0\in{\cal D}_0$ and $(n_k)_k$  a sequence going to $+\infty$ such that $(S^{n_k}(x_0))_k$ 
converges to $x$; we may always assume that $(\TT^{n_k}u_0)_k$ converges to 
$\psi_1$. Consider a sequence $(p_k)_k$ going to $+\infty$ such that $\TT(0,t+n_k+p_k)u_0$ converges uniformly to $u(t)$.\par
An application of Lemma \ref{decroit} tells
us that, for all $\sigma\in\RR$, the function
 $$
 s>0\mapsto \TT(0,s+\sigma) u_0(S^s(x_0))-\phi(s,S^s(x_0))
 $$
 is non-increasing, hence has a
finite limit $l(\sigma)$ as $s\to +\infty$; in particular
 for every $\sigma>0$ the function
$$
s\mapsto\TT(0,\sigma+s+t)u_0(S^{s+t}(x_0))-\phi(s+t,S^{s+t}(x_0))
$$
has the same 
limit $l(\sigma)$ for all $t\in\RR$; we set $l_k=l(p_k)$. We may assume, up to the extraction of a subsequence,
that the sequence $(l_k)_k$ has a limit.\par
On the other hand the weak contraction property implies, for all $t\in [-t_0,t_0]$:
$$
\Vert \TT(0,t+p_k+n_k)u_0-\TT(0,t+p_k)\psi_1\Vert_\infty
\leq \Vert \TT(0,t+n_k)u_0-\TT(0,t)\psi_1\Vert_\infty;
$$
hence $(\TT(0,t+n_k+p_k)u_0-\TT(0,t+p_k)\psi_1)_k$ converges to 0 uniformly on $[-t_0,t_0]
\times\T$. Therefore we have, for all $k$: 
$$\TT(0,p_k+t)\psi_1(S^t(x))-\phi(t,S^t(x))=l_k;
$$ 
letting 
$k\to +\infty$ implies: $u(t,S^t(x))-\phi(t,S^t(x))=l$. We conclude 
by saying that $t_0$ is arbitrary. 
\finproof
It follows that the curve $S^t(x)$ is calibrated by 
$u$,
hence, by Theorem \ref{threg} we obtain:

\begin{cor}
Let $u(t,x)$ be an $\omega$-limit viscosity solution of  (\ref{HJ}).
Then the derivatives $u_x$ and $\phi_x$
exist on $\M$, and they are equal,
$u_x(t,x)=\phi_x(t,x)$
for all $(t,x)\in \M$.
\end{cor}

Before we continue, let us give an important remark.
All the objects constructed in this section, the 
sets $\D$ and $\M$ and the mappings $S^t$,
depend on the periodic solution $\phi$ that 
was chosen in the beginning.
Let us note $\D_{\phi}$ and $\M_{\phi}$ and 
$S^t_{\phi}$ in order to emphasise this dependence.
If $\psi$ is another $1$-periodic viscosity solution,
then we see from Corollary \ref{const} that
the orbits of $\M_{\psi}$ are calibrated by $\phi$.
It follows that the set 
$$
\A=\bigcap _{\phi}\D_{\phi},
$$
is not empty,
where the intersection is taken on the set
of $1$-periodic viscosity solutions.
This set is usually called the Aubry set.
The mappings $S^t_{\phi|\A_0}$
do not depend on $\phi$,
and for all $\phi$, we have 
$\M_{\phi}\subset \A$.

\subsection{Uniqueness set}
We mean a set such that two global solutions of
(\ref{HJ}) coinciding on this set coincide everywhere. 
We formulate the results in the general, non-autonomous, setting.

\begin{prop}
\label{unique}
Let  $u(t,x):\RR\times \T\lto\RR$
be a global and bounded viscosity solution of (\ref{HJ}).
such that the functions $u(t,.),t\in \RR$ are equicontinuous and 
let $\phi(t,x):\RR\times \T\lto\RR$ 
be a 1-periodic viscosity solution of  (\ref{HJ}).
If $u=\phi$ on $\M$,
then $u=\phi$.
\end{prop}
Note that the condition of equicontinuity is satisfied
by $\omega$-limit solutions.

\noindent{\sc Proof.}
Let us fix a point $(s,q)\in \RR\times\T$.
There exists a curve $x(t):( -\infty,s]\lto \T$
which is calibrated by $\phi$ and satisfies
$x(s)=q$.
In view of Lemma \ref{decroit},
the function 
$t\longmapsto u(t,x(t))-\phi(t,x(t))$
is non-increasing hence it has a limit as $t\lto -\infty$.
Since $\gamma$ is calibrated by $\phi$,  
there exists an increasing sequence 
$n_k$ such that $x(-n_k)$ has a limit $x\in \M_0$.
We obtain:
$$
\lim _{t\lto -\infty} u(t,x(t))-\phi(t,x(t))
=\lim _{k\lto \infty} u(-n_k,x(-n_k))-\phi(-n_k,x(-n_k))$$
$$
=\lim u(-n_k,x)-\phi(-n_k,x)=0
$$
where we have
used the equicontinuity of 
the functions $u(t,.)$
and the fact that $u=\phi$ on $\M$.
It follows that 
$
u(t,x(t))-\phi(t,x(t))\leq 0
$
for all $t\leq s$
hence 
$
u(s,q)\leq \phi(s,q).
$
The reversed inequality is obtained in exactly the same
way, with the help
of Lemma \ref{lemuniq} below,  using a curve
 $\gamma(t):( -\infty,s]\lto \T$ calibrated by $u$
and ending at $q$.

\begin{lem}\label{lemuniq}
Let  $u(t,x):\RR\times \T\lto\RR$
be a global and bounded viscosity solution of (\ref{HJ})
and $\gamma(t):]-\infty,s]\lto \T$
be a curve calibrated by $u$.
Then there exists an increasing sequence $n_k$
of integers such that
$\gamma(-n_k)\lto x\in \M_0$.
\end{lem}
\noindent{\sc Proof.}
Let us consider a periodic viscosity solution $\phi$.
In view of Lemma \ref{decroit},
the function 
$t\longmapsto u(t,\gamma(t))-\phi(t,\gamma(t))$
is non-increasing, and bounded.
Hence this function has a limit as $t\lto -\infty$.
Let us choose an increasing  sequence $t_k$ of integers
such that the curves 
$\gamma(t-t_k)$
are converging uniformly on compact sets
to a limit $\gamma_{\infty}:\RR\lto\T$.
The following calculations show that 
this curve is calibrated by $\phi$:
$$
\phi(t',\gamma_{\infty}(t'))-\phi(t,\gamma_{\infty}(t))
=
\lim\Big(
\phi(t',\gamma(t'-t_k))-\phi(t,\gamma(t-t_k))
\Big)
$$
$$
=\lim\Big(
\phi(t'-t_k,\gamma(t'-t_k))-\phi(t-t_k,\gamma(t-t_k))\Big)
$$
$$
=\lim\Big(
u(t'-t_k,\gamma(t'-t_k))-u(t-t_k,\gamma(t-t_k))
\Big)
$$
$$
=\lim \int_{t}^{t'}L(\sigma,\gamma(\sigma-t_k),
\dot\gamma(\sigma-t_k)) d\sigma
=\int_{t}^{t'}L(\sigma,\gamma_{\infty}(\sigma),
\dot\gamma_{\infty}(\sigma)) d\sigma.
$$
As a consequence, the curve $t\longmapsto (t,\gamma_{\infty}(t)) $
is asymptotic to $\M$,
so that there exists an increasing sequences $m_k$
of integers such that
that $\gamma_{\infty}(m_k)\lto x\in \M_0$.
Possibly taking a subsequence of $t_k$, we obtain
that $\gamma(m_k-t_k)\lto x\in \M$
and
that 
$n_k=t_k-m_k$
is increasing.
\finproof

Proposition \ref{unique}
is of course very useful to obtain uniqueness result
for 1-periodic viscosity solutions.
As an example, we obtain that, if there exists a 
periodic solution $\phi$ such that
$S^1_{\M_0}$ has a dense orbit, then
all the other $1$-periodic solutions are of the form
$c+\phi$.
It is a classical result from one-dimensional dynamics that 
this holds on the circle if the rotation
number is irrational, see Section
\ref{circle}.
This yields the uniqueness result of E and Sobolevskii
(\cite{E} and \cite{S}).

In the autonomous case, where the Hamiltonian $H$
does not depend on the variable $t$, the above remarks 
imply that any $\omega$-limit viscosity solution $u$
is independent of $t$ on $\M=\RR\times \M_0$.
One can then conclude that the solution
$u$ is independent of $t$.

The time-periodic case however is more complicated,
and we are not able to give a description
of $\omega$-limit orbits without using 
some specific features of the low dimension.
This will be done in the next section.

\Section{Rotation number and convergence}
\label{circle}
In this section, we shall take advantage of the low dimension.
More precisely, we shall make use of Poincar\'e theory
of homeomorphisms of the circle, see \cite{KH} for example.
We have constructed in the previous section
a closed  subset $\A$ of $\RR\times \T$, 
which is the disjoint union of graphs of calibrated curves.
These calibrated curves will be called the orbits
of $\A$ from now on.
It will be useful to consider the standard projection
$\pi:\RR\times\RR \lto \RR\times \T$,
and the subset $\bar \A=\pi^{-1}(\A).$
Since $\A$ is the disjoint  union of the graphs of its orbits 
$\gamma:\RR\lto \T$, $\bar \A$ is the disjoint
union of graphs of continuous curves 
$\bar \gamma:\RR\lto \RR$.
These curves $\bar \gamma$ are called orbits of $\bar \A$,
they satisfy $\pi\circ \bar \gamma=\gamma$.
By straightforward extension of Poincar\'e Theory
of homeomorphisms of the circle,
we obtain the following result:

The limit 
$
\rho =\di\lim_{t\lto\infty} \bar \gamma(t)/t
$
exists and does not depend on the orbits $\bar \gamma$.
The number $\rho$ is called the rotation number.

\medskip
\noindent{\sc Proof of the theorem.} 
Let $u(t,x)$ be an $\omega$-limit
viscosity solution.
Let $\phi(t,x)=\liminf u(t+n,x)$
be as in Proposition \ref{liminf}.
We will prove that $u=\phi$ in the two cases 
$\rho=0$ and $\rho$ irrational.
The general case of a rational rotation number
$\rho=a/b$, $a\neq 0$  can be reduced to the
case $\rho =0$ by considering the Hamiltonian
$$
\tilde H(t,x,p)=aH(at,bx-at,\frac{p}{b})
$$
and noticing that the function $u(t,x)$
is a solution of the equation (\ref{HJ})
with Hamiltonian $H$ if and only if
the function $\tilde u(t,x)=u(at,bx-at)$
is a solution of the equation (\ref{HJ})
with Hamiltonian $\tilde H$,
and that the rotation number associated 
to this second equation is $\rho=0$.

In view of Proposition \ref{unique}, it is enough
to prove equality on $\M_{\phi}$
(that will be denoted $\M$ from now on).
Let us note $d=u-\phi$ 
and recall from the previous discussions
the main properties of this function $d$.

\medskip
\noindent $\bullet$ The derivative $d_x$ exists on $\M$, and
$d_x(t,x)=0$ for all $(t,x)\in \M$.

\noindent $\bullet$ The function $d$ is constant on the graph of orbits of $\M$.

\noindent $\bullet$ The functions $d(t,.), t\in  \RR$ are equilipschitz.

\noindent $\bullet$ We have 
$\liminf d(n,x)=0$
for each $x\in \T$.

\medskip
\noindent{\bf Case 1.} $\rho=0$.
In this case, the orbits $\gamma$ of $\A$ are of two kinds:\\
(I). Orbits of period $1$.\\
(II). Heteroclinic orbits connecting orbits of type (I).
This means that if $\gamma(t)$ is an orbit of type (II),
there exist orbits $\gamma^{+}$ and $\gamma^-$
of type (I) such that
$\di\lim_{t \lto \infty}(\gamma(t)-\gamma^{+}(t))=0$
and
$\di\lim_{t \lto -\infty}(\gamma(t)-\gamma^{-}(t))=0$.

It follows that the $\omega$-limit $\M$ is the union
of the graphs of the orbits of $\A$ of type (I).
Let $\gamma$ be such an orbit.
The function $d$ is constant on the graph of $\gamma$.
Since in addition we have 
$0=\liminf d(n,\gamma(0))=\liminf d(n,\gamma(n))$,
this constant is zero.
As a consequence, $d=0$ on $\M$.

\medskip
\noindent {\bf Case 2.} $\rho$ is irrational. 
In this case, there is no periodic orbit in $\A$,
hence in $\M$.
We have to describe the complement of $\M$:
\begin{lem}
Each connected component $\bar U$
 of the  complement 
of $\bar \M$ in $\RR\times\RR$ is of the form
$$
\big\{
(t,x)  \hbox{ such that }
\bar \gamma^-(t)\leq x \leq \bar  \gamma^+(t)
\big\},$$
where $\bar \gamma^{\pm}(t)$ are two orbits of $\bar \M$
which satisfy $\di\lim_{t\to\pm\infty} |\bar \gamma^+(t)-\bar \gamma^-(t)|=0$.
\end{lem}
\noindent{\sc Proof.}
It is quite clear that each connected    component $\bar U$
 of the  complement 
of $\bar \M$ in $\RR\times\RR$ is of the form
$
\bar \gamma^-(t)\leq x \leq \bar  \gamma^+(t)
.$
In order to prove that
 $\di\lim_{t\to\pm\infty} |\bar \gamma^+(t)- \bar\gamma^-(t)|=0$,
let us consider the interval 
$\bar I_k=]\bar \gamma^-(k),\bar \gamma^+(k)[\subset \RR$.
We claim that the intervals $I_k=\pi(\bar I_k)$ are all disjoint
in $\T$.
Indeed, recalling that $\M_0$ is the set of points 
$x\in \T$ such that $(0,x)\in \M$, we see that each of the
intervals $I_k$ is a connected component in $\T$ of the complement
of $\M_0$. 
Now suppose that $I_k$ and $I_l$ have nonempty intersection.
then the boundaries of $I_k  \cap  I_l$ are points
of $\M_0$, (since the boundaries of $I_k$ and $I_l$ are)
so they are contained neither in $I_k$ nor in $I_l$,
which is possible only if $I_k= I_l$.
Now assume that $I_k=I_l$ with $k < l$.
This implies that $\gamma^+(l)= \gamma^+(k)$,
hence $\gamma^+(t) $ is periodic, which is a contradiction.
So the intervals $I_k$ are all disjoint, hence the sum
$\di\sum_k(\bar \gamma^+(k)-\bar \gamma^-(k))$ 
of their lengths is finite; hence the result.
\finproof

The proof that $d|_{\M}=0$ is similar to the proof
of Proposition 6.5. in \cite{Be}.
Let us set $\bar d=d\circ \pi$.
Let $\bar U$ be a connected component of  the
complement of $\bar \M$, and let $\bar \gamma^{\pm}$
be its boundary curves.
Clearly, $\bar d(t,\bar \gamma^+(t))$ and 
 $\bar d(t,\bar \gamma^-(t))$ are constants
$c^+$ and $c^-$.
In addition, since $\di\lim_{t\to\pm\infty} ( \bar \gamma^+(t) -\bar \gamma^-(t)=0$
and since the functions $\bar d(t,.), t\in \RR$ are equilipschitz,
we have  $c^+=c^-$.
Since this holds for every connected component $\bar U$
of the complement, there exists a continuous function 
$D:\RR\times\RR\lto \RR$ which coincides with
$\bar d$ on $\bar \M$ and is constant on the complement
of $\bar \M$.
Since $d_x=0$ on $\M$,
the derivative $D_x$ exists and vanishes at each point.
It follows that the function $D$ is constant,
hence $\bar d$ is constant on $\bar \M$, hence
$d$ is constant on $\M$.
Recalling that $\liminf d(n,x)=0$ for $x\in \M_0$,
we obtain $d=0$ on $\M$.
\finproof
\textbf{Remark.}
The discussions above in case 2 also yields that, 
if $\phi$ and $\psi$ are two periodic visicosity solutions,
then $d=\psi-\phi$ is constant on $\M$, 
hence $\psi-d $ and $\phi$ are two periodic viscosity solutions
which are equal on $\M$, hence equal by Proposition \ref{unique}.
This provides a proof of the uniqueness result of 
E and Sobolevskii.
Note however that this result more simply follows
from the fact that all the orbits of $S^1_{\M_0}$
are dense in $\M_0$.

\end{document}